\documentclass{amsart}
\usepackage{amsthm}
\usepackage{amssymb}

\theoremstyle{plain}
\newtheorem{thm}{Theorem}[section]
\newtheorem{prop}[thm]{Proposition}
\newtheorem{cor}[thm]{Corollary}
\newtheorem{lem}[thm]{Lemma}

\theoremstyle{remark}
\newtheorem*{remark}{Remark}
\newtheorem*{exmp}{Example}

\numberwithin{equation}{section}

\date{Probability Theory and Related Fields, Received: 5 August 2001, 
Revised: 28 March 2007, Published online: 27 April 2007}

\title[\null]{AR and MA representation of 
partial autocorrelation functions, 
with applications}
\author[\null]{Akihiko INOUE}

\address{Department of Mathematics \\
Faculty of Science \\
Hokkaido University \\
Sapporo 060-0810 \\
Japan}
\email{inoue@math.sci.hokudai.ac.jp}


\begin{document}

\begin{abstract}
We prove a representation of the partial autocorrelation function (PACF), 
or the Verblunsky coefficients, of a stationary process 
in terms of the AR and MA coefficients. We apply it to 
show the asymptotic behaviour of the PACF. 
We also propose a new definition of short and long memory in terms of 
the PACF.
\end{abstract}

\maketitle

%%%%%%%%%%%%%%%%%%%%      Section 1.    %%%%%%%%%%%%%%%%%%
\section{Introduction}\label{sec:1}

Let $\{X_n: n\in\mathbf{Z}\}$ be a real, zero-mean, weakly stationary 
process, defined on a probability space $(\Omega,\mathcal{F},P)$ 
with spectral measure not of finite support, which we 
shall simply call a {\it stationary process\/}. 
We write $\{\gamma_n: n\in\mathbf{Z}\}$ for 
the autocovariance function of $\{X_n\}$: 
$\gamma_n:=E[X_nX_0]$ for $n\in\mathbf{Z}$. 
For $\{X_n\}$, we have another sequence 
$\{\alpha_n\}_{n=0}^{\infty}$ 
called the {\it partial autocorrelation function\/} (PACF), 
where 
$\alpha_0:=\gamma_0$, $\alpha_1:=\gamma_1/\gamma_0$, and for $n\ge 2$, 
$\alpha_n$ is the correlation coefficient of 
the two residuals obtained from $X_0$ and $X_n$ by regressing on the 
intermediate values $X_1, \ldots\ X_{n-1}$ (see \S 2 below). 

The autocovariance function $\{\gamma_n\}$ is positive definite, and the 
inequalities that this positive definiteness imposes may be inconvenient 
in some contexts. 
By contrast, the PACF $\{\alpha_n\}_{n=0}^{\infty}$ 
gives an {\it unrestricted parametrization\/}, in that the only inequalities 
restricting the $\alpha_n$ are the obvious ones implied by their being 
correlation coefficients, i.e., $\alpha_n \in [-1,1]$, 
or $(-1,1)$ in the non-degenerate case relevant here. 
This result is due to Barndorff--Nielsen and Schou \cite{BS}, 
Ramsey \cite{Ra} in the time-series context.  See also D\'egerine \cite{De}, 
and for extensions to the non-stationary case, 
D\'egerine and Lambert--Lacroix \cite{DL}. 
However, in the context of mathematical analysis --- specifically, 
the theory of orthogonal polynomials on the unit circle (OPUC) --- the result 
dates back to 1935-6 to work of Verblunsky \cite{V1, V2}, 
where the PACF appears as the sequence of {\it Verblunsky coefficients\/}. 
For a survey of OPUC, see Simon \cite{Si2}, and for a textbook treatment, 
\cite{Si3} (analytic theory), \cite{Si4} (spectral theory).  One of our main purposes here is to emphasize the importance of the PACF: ``This knows everything", a remark we owe to Yukio Kasahara.

The question thus arises of a `dictionary', allowing one to pass between 
statements on the covariance $\{\gamma_n: n\in\mathbf{Z}\}$, 
or the spectral measure $\mu$ defined by 
$\gamma_n=\int_{-\pi}^{\pi}e^{ik\theta}\mu(d\theta)$, 
and the PACF $\{\alpha_n\}_{n=0}^{\infty}$; 
see the last paragraph of 
\cite{Si3}, page 3, where this is stated as perhaps ``the 
most central'' question in OPUC. 
A prototype of such a result is {\it Baxter's theorem\/} 
(see \cite{Ba1,Ba2}; see also Theorem \ref{thm:4.1} below and 
\cite[Chapter 5]{Si3}). 
Such a link exists, in the shape of the Levinson 
(or Szeg\"o--Levinson--Durbin) 
algorithm, due to Szeg\"o \cite{Sz} in 1939 (Szeg\"o recursion), 
Levinson \cite{L} in 1947, 
Durbin \cite{Du} in 1960; 
for a textbook 
account, see Pourahmadi \cite{P}, \S 7.2. 
For the history of Szeg\"o recursion (Theorem 1.5.2 of \cite{Si3}), 
see \cite{Si3}, page 69. 
But while very useful numerically, 
the Levinson algorithm is less suitable for theoretical studies, such as one 
of the questions that motivates us here --- the behaviour of the PACF for 
large lags (that is, large $n$).

In this paper, to fill the gap above at least partially, 
we introduce a representation of the PACF which is given only in terms of 
another sequence $\{\beta_n\}_{n=0}^{\infty}$ defined by
\begin{equation*}
\beta_n,\ \mbox{or}\ 
\beta(n),\ :=\sum_{v=0}^{\infty} c_v a_{v+n}\qquad(n=0,1,\dots)
\tag{$\ast$}
\end{equation*}
(see Inoue and Kasahara \cite[page 8]{IK1}, \cite[(2.23)]{IK2}), 
where the two sequences $\{c_n\}_{n=0}^{\infty}$ and 
$\{a_n\}_{n=0}^{\infty}$ have statistical interpretations as the coefficients 
of the MA$(\infty)$ and AR$(\infty)$ 
representations of the process, respectively 
(see (MA) and (AR) below). 
They also have analytic interpretations as the coefficients in the Maclaurin 
expansions of the {\it Szeg\"o function\/} $D(z)$ occurring in the 
theory of OPUC and its associate $-1/D(z)$, 
both of which are {\it outer functions} in Beurling's sense 
(see \S 2 below). 
For background, see \cite{Si3}, 
or --- from a statistical point of view --- 
Grenander and Szeg\"o \cite{GS}, Rozanov \cite{Ro}, 
Ibragimov and Rozanov \cite{IR}.

We are particularly interested in the asymptotics of $\alpha_n$ as 
$n \to \infty$, and the representation of $\alpha_n$ is useful 
in investigating them since the representation enables us to study 
$\alpha_n$ directly via $\beta_n$ or $a_n$ and $c_n$. 
In a number of the specific examples of processes 
with {\it long memory\/} we treat, 
here and in \cite{I2}, \cite{I3}, \cite{IK1}, we observe behaviour of the form
\begin{equation*}
\alpha_n \sim d/n \qquad (n \to \infty)
\tag{$d/n$}
\end{equation*}
(by the representation of $\alpha_n$, 
we are able to improve the estimate of this type in \cite{I2,I3,IK1}, 
where only $\vert \alpha_n\vert$ was considered --- 
we were unable to determine its sign). 
In ($d/n$), and throughout the paper, 
$a_n \sim b_n$ as $n\to\infty$ means 
$\lim_{n\to\infty}a_n/b_n=1$. 
On the one hand, ($d/n$) seems very special behaviour, 
if we begin by specifying our model via the PACF, since by 
Verblunsky's theorem any value in $(-1,1)$ can be taken by any $\alpha_n$. 
On the other hand, it is more usual in practice to specify our model 
with long memory by other 
means, such as the AR and MA coefficients, and here such 
behaviour seems to be typical of the (quite broad) classes of example 
where one can carry out the computations and obtain an explicit asymptotic 
expression for $\alpha_n$, such as the fractional ARIMA (FARIMA) models 
studied in \cite{I3, IK1} and Theorem \ref{thm:2.5} below. 
In this connection, we note that in recent work of Simon \cite{Si1}, 
the case ($d/n$) is described as `a prototypical example'. 
See also \cite{Si3}, \cite{Si4} where \cite{Si1} is developed further, 
and the papers by Golinskii and Ibragimov \cite{GI}, 
Damanik and Killip \cite{DK} that inspired it. We note also that a 
generalization of ($d/n$), in which `asymptotic to' is replaced by `of the 
same order of magnitude as', is obtained 
as the conclusion in work of Ibragimov and Solev \cite{IS}, 
under conditions on the spectral density.

We denote by $H$ the real Hilbert space spanned by $\{X_k:k\in\mathbf{Z}\}$ 
in $L^2(\Omega,\mathcal{F},P)$. The norm of $H$ is given by 
$\Vert Y\Vert:=E[Y^2]^{1/2}$. For $n, m\in\mathbf{N}$ with $n\le m$, 
we denote by 
$H_{[n,m]}$, $H_{(-\infty,m]}$ and $H_{[m,\infty)}$ 
the subspaces of $H$ spanned by $\{X_{n},\dots,X_{m}\}$, 
$\{X_k: k\le m\}$, and $\{X_k: k\ge n\}$, respectively. 
The proof of the representation of the PACF is based on the approach 
introduced in \cite{I2} which combines von Neumann's alternating 
projection theorem (see \cite[Theorem 9.20]{P}) and the 
following {\it intersection of past and future\/} 
property of $\{X_n\}$:
\begin{equation*}
H_{(-\infty,n]}\cap H_{[1,\infty)}=H_{[1,n]}\qquad (n=1,2,\dots).
\tag{\textrm{IPF}}
\end{equation*}
This approach is also useful in continuous-time models; 
see Inoue and Nakano \cite{IN} and the references therein. 
A useful sufficient condition for (IPF) is that $\{X_n\}$ 
is purely non-deterministic (PND) (see \S 2 below) and has 
spectral density $\Delta(\cdot)$ such that 
$\int_{-\pi}^{\pi}d\theta/\Delta(\theta)<\infty$ 
(see \cite{I2}, Theorem 3.1), which itself is a discrete-time 
version of the Seghier--Dym theorem \cite{S}, \cite{Dy}. 
This theorem itself originates in work of Levinson and 
McKean \cite{LM}. 
Naturally, (IPF) is closely related to 
the property
\begin{equation*}
H_{(-\infty,0]}\cap H_{[1,\infty)}=\{0\}
\tag{\textrm{CND}}
\end{equation*}
called {\it complete non-determinism\/}; 
see Bloomfield et al.\ \cite{BJH}. 
In fact, a stationary process is completely non-deterministic 
if and only both (PND) and (IPF) are satisfied (see 
\cite{IK2}, Theorem 2.3). In particular, if $\{X_n\}$ 
is PND and has 
spectral density $\Delta(\cdot)$ such that 
$\int_{-\pi}^{\pi}d\theta/\Delta(\theta)<\infty$, then it is CND.

In \S 2, we state the main results, including the 
representation of the PACF and its asymptotic behaviour 
of the type ($d/n$). 
Sections 3--5 are devoted to their proofs. 
In \S 6, we give a further application of the representation, 
that is, the asymptotics for the PACF of processes with 
regularly varying covariance functions. 
In \S 7, we close the paper with results 
for ARMA processes.

%%%%%%%%%%%%%%%%%%%%      Section 2.    %%%%%%%%%%%%%%%%%%
\section{Main Results}\label{sec:2}

As stated in \S 1, 
let $H$ be the real Hilbert space spanned by $\{X_k:k\in\mathbf{Z}\}$ 
in $L^2(\Omega,\mathcal{F},P)$, which has inner product 
$(Y_1,Y_2):=E[Y_1Y_2]$ and norm $\Vert Y\Vert:=(Y,Y)^{1/2}$. 
Also, for an interval $I\subset \mathbf{Z}$, we write $H_I$ for 
the closed subspace of $H$ spanned by $\{X_k: k\in I\}$ 
and $H_{I}^{\bot}$ for the orthogonal complement of $H_I$ in $H$. 
Let $P_I$ and $P^{\bot}_I$ be the orthogonal projection operators of $H$ 
onto $H_I$ and $H_{I}^{\bot}$, respectively. 
Thus $P^{\bot}_IY=Y-P^{\bot}Y$ for $Y\in H$. 
The projection $P_IY$ stands for 
the best linear predictor of $Y$ based on the observations 
$\{X_k: i\in I\}$, and $P_I^{\bot}Y$ for its prediction error.

The partial autocorrelation function (PACF) $\{\alpha_n\}_{n=0}^{\infty}$ of 
$\{X_n\}$ is defined by
\begin{equation*}
\alpha_0:=\gamma_0,\qquad 
\alpha_n:=U_n/V_n\qquad(n=1,2,\dots),
\tag{PACF}
\end{equation*}
where $U_1:=(X_1, X_0)$, $V_1:=\Vert X_1\Vert^2$ and
\begin{align*}
U_n:&=
(P^{\bot}_{[1,n-1]}X_n, P^{\bot}_{[1,n-1]}X_0)
\qquad (n=2,3,\dots),\\
V_n:&=
\Vert P^{\bot}_{[1,n-1]}X_n\Vert^2\qquad (n=2,3,\dots)
\end{align*}
(cf.\ Brockwell and Davis \cite[\S 3.4 and \S 5.2]{BD}). 
We have $V_n>0$ for $n\ge 1$ since we have assumed that 
the spectral measure $\mu$ of $\{X_n\}$ has infinite support; 
$X_1,X_2,X_3,\dots$ are linearly independent since 
no nonzero trigonometric polynomial vanishes $\mu$ a.e. 
Also, since $\Vert X_0\Vert=\Vert X_1\Vert$ and 
$\Vert P^{\bot}_{[1,n-1]}X_0\Vert=\Vert P^{\bot}_{[1,n-1]}X_n\Vert$ 
for $n\ge 2$, 
$\alpha_n$ is actually the correlation coefficient between the residuals 
$P^{\bot}_{[1,n-1]}X_0$ and $P^{\bot}_{[1,n-1]}X_n$ 
(resp.\ $X_0$ and $X_n$) for $n\ge 2$ (resp.\ $n=1$).

A stationary process $\{X_n\}$ is said to be 
{\it purely nondeterministic} (PND) if
\begin{equation*}
\bigcap\nolimits_{n=-\infty}^{\infty}H_{(-\infty,n]}=\{0\}
\tag{PND}
\end{equation*}
or, equivalently, there exists a positive even and 
integrable function $\Delta(\cdot)$ on $(-\pi,\pi)$ such that 
$$
\gamma_n=\int_{-\pi}^{\pi}e^{in\theta}\Delta(\theta)d\theta\quad
(n\in\mathbf{Z}),\qquad 
\int_{-\pi}^{\pi}\vert\log \Delta(\theta)\vert d\theta<\infty
$$
(see \cite[\S 5.7]{BD}, \cite[Chapter II]{Ro} and 
\cite[Chapter 10]{GS}). 
We call $\Delta(\cdot)$ the {\it spectral density\/} of $\{X_n\}$. 
Using $\Delta(\cdot)$, we define the 
{\it Szeg\"o function} $D(\cdot)$ by
\begin{equation*}
D(z):=\sqrt{2\pi}\mbox{exp}\left\{\frac{1}{4\pi}\int_{-\pi}^{\pi}
\frac{e^{i\theta}+z}{e^{i\theta}-z}
\log \Delta(\theta)d\theta\right\}\qquad(z\in\mathbf C,\ \vert z\vert<1).
\end{equation*}
The function $D(z)$ is 
an outer function in the Hardy space $H^{2}$ 
of class 2 over the unit disk $\vert z\vert<1$. 
Using $D(\cdot)$, we define the MA coefficients $c_n$ 
by
\begin{equation*}
D(z)=\sum_{n=0}^{\infty}c_nz^n\qquad(\vert z\vert<1),
\end{equation*}
and the AR coefficients $a_n$ by
\begin{equation*}
-\frac{1}{D(z)}=\sum_{n=0}^{\infty}a_nz^n\qquad(\vert z\vert<1)
\end{equation*}
(see \cite[\S 4]{I2} and \cite[\S 2.2]{IK2} for background). 
Both $\{c_n\}$ and $\{a_n\}$ are 
real sequences, and $\{c_n\}$ is in $l^2$. 
The coefficients $c_n$ and $a_n$ are actually those that appear in the 
following MA$(\infty)$ and AR$(\infty)$ representations, respectively, of 
$\{X_n\}$ (under suitable condition such as $\{a_n\}\in l^1$ for the 
latter):
\begin{align*}
&X_n=\sum_{j=-\infty}^{n}c_{n-j}\xi_j\qquad (n\in\mathbf{Z}),
\tag{MA}\\
&\sum_{j=-\infty}^{n}a_{n-j}X_{j}+\xi_n=0\qquad (n\in\mathbf{Z}),
\tag{AR}
\end{align*}
where $\{\xi_k\}$ is the innovation process defined by 
$\xi_k=\epsilon_k/\Vert \epsilon_k\Vert$ 
with $\epsilon_k$ the prediction error when we predict $X_k$ from 
the whole past $\{X_m: m\le k-1\}$, 
i.e., $\epsilon_k=P^{\bot}_{(-\infty,k-1]}X_k$ (cf.~\cite[\S 2]{IK2}). 
From (MA), we have the following equality:
\begin{equation}
\gamma_n=\sum_{v=0}^{\infty}c_vc_{\vert n\vert+v}\qquad(n\in\mathbf{Z}).
\label{eq:2.1}
\end{equation}

We wish to derive a representation 
of $\alpha_n$ which is given only in terms of $\beta(\cdot)$ defined by 
($\ast$). 
For this purpose, we consider the following two conditions 
(BC) and (O($1/n$)):

\ 

\noindent
(BC)\quad The process $\{X_n\}$ has summable autocovariance function 
$\{\gamma_n\}$, i.e., 
$\sum_{-\infty}^{\infty}\vert \gamma_n\vert<\infty$, and {\it positive\/} 
spectral density $\Delta(\cdot)$, i.e., 
$\min_{\theta\in [-\pi, \pi]} \Delta(\theta)>0$.

\ 

\noindent
(O($1/n$))\quad $\{X_n\}$ is PND, and satisfies both 
$\sum_{0}^{\infty}\vert a_n\vert<\infty$ 
and
\begin{equation}
\sum_{v=0}^{\infty}\vert c_va_{n+v}\vert=O(1/n)\qquad(n\to\infty).
\label{eq:2.2}
\end{equation}

\

Notice that if $\{\gamma_n\}\in l^1$, 
then $\{X_n\}$ has continuous spectral density 
$\Delta(\theta)=(2\pi)^{-1}\sum_{-\infty}^{\infty}\gamma_ke^{-ik\theta}$. 
We will see (Theorem \ref{thm:4.1}) that (BC) holds 
if and only if 
$\{X_n\}$ is PND, $\{a_n\}\in l^1$ and $\{c_n\}\in l^1$. 
We also see that (O($1/n$)) holds for many interesting processes including 
the FARIMA$(p,d,q)$ processes with $0<d<1/2$ 
which we consider below. 
The condition (L$(d,\ell)$) below imples (O($1/n$)) 
(see Proposition \ref{prop:5.1} below).
ARMA processes satisfy both (BC) and (O($1/n$)).

\ 

\noindent {\it Remarks.}\quad 
1.\ 
Condition (BC) requires that the process have {\it summable autocovariance} 
and {\it positive spectral density}.  It could thus be denoted (SP).  We call it (BC) instead 
to emphasize its role as {\it Baxter's condition}. 
{\it Baxter's theorem} (\cite{Ba1,Ba2}; see also Theorem \ref{thm:4.1} below) 
gives the equivalence of (BC) with summability of 
the PACF (i.e., $\{\alpha_n\} \in l^1$) 
subject only to the very weak condition that the 
spectral measure $\mu$ has infinite support.  
See Chapter 5 of [Si3], where Baxter's theorem is discussed in detail and 
proved.

\noindent 2.\ 
Several different definitions of {\it short memory} (or its negation, {\it long memory}) are in current use.  
See for example Section 2 of the survey paper Baillie \cite{Bai} 
for details and references.  The most standard definition is that $\{X_n\}$ has long memory (resp.\ short memory) 
if $\sum_{k=-\infty}^{\infty}\vert\gamma_k\vert=\infty$ 
(resp.\ $<\infty$); see Beran \cite[page 6]{Be}, 
\cite[\S 13.2]{BD}.  The fact that Baxter's theorem is so powerful and useful suggests the possibility of using Baxter's condition to define short memory in a new sense: call the process {\it short memory} if Baxter's condition holds, {\it long memory} otherwise. 

\noindent 3.\ The difference between these approaches to short and long memory is well illustrated by the fractional ARIMA (or FARIMA) processes (see below for definitions), studied in e.g. \cite{I3,IK1} and 
Theorems \ref{thm:2.4} and \ref{thm:2.5} below.  
The two main cases $d \in (-1/2,0)$ 
and $d \in (0, 1/2)$ behave in the {\it same} way from the point of 
view of asymptotics of PACF ($d/n$ for each; see Theorem \ref{thm:2.5} 
below) and 
prediction error 
($d^2/n$ for each; see (\ref{eq:2.24}) below) ---  
but in {\it different} ways from the point of 
view of summability of the autocovariance function.  
Our contention is that the PACF 
$\{\alpha_n\}$, and/or the prediction error $\{\delta(n)\}$, 
are more informative about the essence of long-range dependence --- 
the rate at which the information in the remote past decays with time --- 
than the autocovariance function 
$\{\gamma_n\}$ which is usually used here.

\noindent 4.\ 
Wu \cite{W}, \S 3 considers stationary processes which 
satisfy $\{\gamma_n\}\notin l^1$ (long-range dependence) but 
have $\lim_{K\to\infty}\sum_{k=-K}^K\gamma(k)$. 
He uses the Zygmund class of slowly varying functions; see 
Bingham et al.\ \cite[\S 1.5.3]{BGT}, Zygmund \cite[V.2]{Z}.

\ 

Under (BC) or (O($1/n$)), we define, for $n=0,1,\dots$,
\begin{align}
d_1(n)&=\beta(n),
\label{eq:2.3}\\
d_2(n)&=\sum_{m_1=0}^{\infty}\beta(m_1+n)\beta(m_1+n),
\label{eq:2.4}
\end{align}
and, for $k=3,4,\dots$,
\begin{equation}
\begin{split}
&d_k(n)
=\sum_{m_{k-1}=0}^{\infty}
\beta(m_{k-1}+n)
\sum_{m_{k-2}=0}^{\infty}\beta(m_{k-1}+m_{k-2}+n) \\
&\qquad\cdots \sum_{m_2=0}^{\infty}\beta(m_3+m_2+n)
\sum_{m_1=0}^{\infty}
\beta(m_2+m_1+n)\beta(m_1+n),
\end{split}
\label{eq:2.5}
\end{equation}
the sums converging absolutely (see Proposition \ref{prop:4.3} below).

Here is the representation of the PACF $\{\alpha_n\}$.

\begin{thm}\label{thm:2.1}
We assume either $(\mathrm{BC})$ or {\rm (O($1/n$))}. 
Then, for $n=1,2,\dots$, 
\begin{align}
U_n&=(c_0)^2\sum\nolimits_{k=1}^{\infty}d_{2k-1}(n),
\label{eq:2.6}\\
V_n&=(c_0)^2\left\{1+\sum\nolimits_{k=1}^{\infty}d_{2k}(n)\right\},
\label{eq:2.7}\\
\alpha_n&=\frac{\sum_{k=1}^{\infty}d_{2k-1}(n)}
{1+\sum_{k=1}^{\infty}d_{2k}(n)},
\label{eq:2.8}
\end{align}
all the sums converging absolutely.
\end{thm}

If $\{X_n\}$ is PND, then $V_n\downarrow 
\Vert P^{\bot}_{(-\infty,-1]}X_0\Vert^2=(c_0)^2$, whence 
$\alpha_n=U_n/V_n\sim (c_0)^{-2}U_n$, 
as $n\to\infty$ (see \cite[\S 2]{I2}). 
Thus an immediate consequence of Theorem \ref{thm:2.1} is 
the next corollary.

\begin{cor}\label{cor:2.2}
We assume either $(\mathrm{BC})$ or {\rm (O($1/n$))}. 
Then, 
\begin{equation}
\alpha_n \sim \sum_{k=1}^{\infty}d_{2k-1}(n)\qquad 
(n\to\infty).
\label{eq:2.9}
\end{equation}
\end{cor}

We turn to the results on the asymptotic behaviour of $\alpha_n$ as 
$n\to\infty$. We write $\mathcal{R}_0$ for the class of 
slowly varying functions at infinity: the class of positive, measurable 
$\ell$, defined on some neighborhood $[A,\infty)$ of infinity, such that
\begin{equation}
\lim_{x\to\infty}\ell(\lambda x)/\ell(x)=1\quad\mbox{for all $\lambda>0$}
\label{eq:2.10}
\end{equation}
(see \cite[Chapter 1]{BGT}). 
For $\ell\in\mathcal{R}_0$ and $d\in (0,1/2)$, 
we consider the following condition as a standard 
one for processes with {\it long memory\/} 
(see \cite[(A2)]{IK2}): 

\ 

\noindent 
(L$(d,\ell)$)\quad $\{X_n\}$ is PND and $\{c_n\}$ and $\{a_n\}$ 
satisfy, respectively,
\begin{align}
&c_n\sim n^{-(1-d)}\ell(n)
\qquad(n\to\infty),
\label{eq:2.11}\\
&a_n\sim n^{-(1+d)}\frac{1}{\ell(n)}
\frac{d\sin(\pi d)}{\pi}
\qquad (n\to\infty).
\label{eq:2.12}
\end{align}

\ 

The condition (L$(d,\ell)$) implies 
\begin{equation}
\gamma_n\sim n^{-(1-2d)}\ell(n)^2B(d,1-2d)\qquad (n\to\infty)
\label{eq:2.13}
\end{equation}
(see \cite[(2.22)]{IK2}), whence 
$\{\gamma_n\}\notin l^1$, where $B(\cdot,\cdot)$ denotes the beta integral, 
i.e.,
\[
B(p,q)=\int_0^1t^{p-1}(1-t)^{q-1}dt=\frac{\Gamma(p)\Gamma(q)}{\Gamma(p+q)}
\qquad (p,q>0).
\]

Here is a result of the type ($d/n$) for processes with 
long memory.

\begin{thm}\label{thm:2.3}
For $\ell\in\mathcal{R}_0$ and $d\in (0,1/2)$, 
$(L(d,\ell))$ implies the asymptotic behaviour $(d/n)$ of the 
PACF.
\end{thm}

For example, the FARIMA$(p,d,q)$ 
model with $0<d<1/2$, which is regarded as a standard parametric model 
with long memory and which we consider below, satisfies 
(L$(d,\ell)$) for some constant function $\ell$ 
(see Kokoszka and Taqqu \cite[Corollary 3.1]{KT}), 
whence, by this theorem, its PACF has the asymptotic behaviour $(d/n)$. 
We see another class satisfying (L$(d,\ell)$) in \S \ref{sec:6}. 
We can write the property ($d/n$) as $d=\lim_{n\to\infty}n\alpha_n$, 
suggesting how to estimate the parameter $d$, which is 
important in a process with long memory. See \cite[\S 5]{IK1} for 
numerical calculation.

We consider the FARIMA. For 
$d\in (-1/2, 1/2)$ and $p, q\in\mathbf{N}\cup\{0\}$, 
a stationary process $\{X_n\}$ is said to be a (causal and invertible) 
fractional 
ARIMA$(p,d,q)$ (or FARIMA$(p,d,q)$) process if it has a spectral 
density $\Delta(\cdot)$ of 
the form
\begin{equation}
\Delta(\theta)
=\frac{1}{2\pi}
\frac{\vert \Theta(e^{i\theta})\vert^2}{\vert \Phi(e^{i\theta})\vert^2}
\vert 1-e^{i\theta}\vert^{-2d}\qquad(-\pi<\theta<\pi),
\label{eq:2.14}
\end{equation}
where $\Phi(z)$ and $\Theta(z)$ are polynomials with real coefficients 
of degrees $p$, $q$, respectively, satisfying the following condition:
\begin{equation}
\begin{split}
&\mbox{$\Phi(z)$ and $\Theta(z)$ have no common zeros, and have no zeros}\\
&\mbox{in the closed unit disk $\{z\in\mathbf{C}: \vert z\vert\le 1\}$.}
\end{split}
\label{eq:2.15}
\end{equation}
The fractional ARIMA model was introduced independently 
by Granger and Joyeux \cite{GJ} and Hosking \cite{Ho}. 
See \cite[\S 2.5]{Be} and \cite[\S 13.2]{BD} for textbook 
treatments and \cite{KT} for a formulation in terms of 
the backward shift operator $B$. 
If $d\in (-1/2, 1/2)\setminus\{0\}$, then the 
MA coefficient $c_n$, AR coefficients $a_n$ and 
autocovariance function $\{\gamma_n\}$ 
of the FARIMA$(p,d,q)$ process 
$\{X_n\}$ satisfy
\begin{align}
&c_n\sim n^{-(1-d)}\frac{K_1}{\Gamma(d)}\qquad(n\to\infty),
\label{eq:2.16}\\
&a_n\sim n^{-(1+d)}\frac{\Gamma(d)}{K_1}\cdot\frac{d\sin(\pi d)}{\pi}
\qquad(n\to\infty),
\label{eq:2.17}\\
&\gamma_n\sim n^{-(1-2d)}
\frac{(K_1)^2\Gamma(1-2d)\sin(\pi d)}{\pi}
\qquad(n\to\infty),
\label{eq:2.18}
\end{align}
where
\begin{equation}
K_1:=\Theta(1)/\Phi(1)
\label{eq:2.19}
\end{equation}
(see \cite[Corollary 3.1]{KT} and \cite[\S 3]{I3}). 
In particular, if $0<d<1/2$, then $\{X_n\}$ has 
long memory. On the other hand, if $d=0$, 
then $\{X_n\}$ reduces to the ordinary ARMA$(p,q)$ process, 
and each of $c_n$, $a_n$ and 
$\gamma_n$ decays exponentially 
fast as $n\to\infty$ (see \cite[Chapter 3]{BD} and \S 7 below).

Theorems \ref{thm:2.1} and \ref{thm:2.3} cover FARIMA with $d\in (0,1/2)$ 
but not the case $d\in (-1/2,0)$. 
However, the FARIMA$(0,d,0)$ with $d\in (-1/2,1/2)$ 
has the PACF given by
\begin{equation}
\alpha_n=\frac{d}{n-d}\qquad(n=1,2,\dots)
\label{eq:2.20}
\end{equation}
(see \cite[Theorem 1(f)]{Ho} as well as \cite[(13.2.10)]{BD}), 
which suggests that ($d/n$) may also hold even if $d\in (-1/2,0)$, 
and this is actually the case as we show below.

For FARIMA with $-1/2<d<0$, we define
\begin{equation*}
\phi_n:=
\begin{cases}
-a_0&\quad (n=0),\\
a_{n-1}-a_{n}& \quad(n=1,2,\dots).
\end{cases}
\end{equation*}
and
\begin{equation*}
\psi_n:=-\sum_{k=n+1}^{\infty}c_k\qquad(n=0,1,\dots).
\end{equation*}
Notice that $\phi_n$ here corresponds to $-\phi_n$ in \cite{IK1}. 
We define $q\in (1/2,1)$ by
\[
q:=1+d.
\]
Then, by \cite[Lemma 4.1]{IK1} and \cite[Corollary 3.1]{KT} 
(see also \cite[Lemma 2.1]{I3}), we have
\begin{gather}
\psi_n\sim n^{-(1-q)}\frac{K_1}{\Gamma(q)}
\qquad(n\to\infty),
\label{eq:2.21}
\\
\phi_n\sim -n^{-(1+q)}\frac{\Gamma(q)}{K_1}\cdot 
\frac{q\sin(\pi q)}{\pi}
\qquad (n\to\infty),
\label{eq:2.22}
\end{gather}
where $K_1$ is as in (\ref{eq:2.19}). 
We define
$$
\beta_-(n):=\sum_{v=0}^{\infty}\psi_v\phi_{v+n+1}\qquad(n=0,1,\dots).
$$
Notice that $\beta_-(n)$ corresponds to $-\beta(n)$ in \cite{IK1}. 
We define, for $n=0,1,\dots$,
\begin{align*}
d_1(n)&=\beta_-(n),
\\
d_2(n)&=\sum_{m_1=0}^{\infty}\beta_-(m_1+n)\beta_-(m_1+n),
\end{align*}
and, for $k=3,4,\dots$,
\begin{equation*}
\begin{split}
&d_k(n)=\sum_{m_{k-1}=0}^{\infty}
\beta_-(m_{k-1}+n)
\sum_{m_{k-2}=0}^{\infty}\beta_-(m_{k-1}+m_{k-2}+n) \\
&\qquad\cdots \sum_{m_2=0}^{\infty}\beta_-(m_3+m_2+n)
\sum_{m_1=0}^{\infty}
\beta_-(m_2+m_1+n)\beta_-(m_1+n),
\end{split}
\end{equation*}
the sums converging absolutely (see \cite[Theorem 3.3]{IK2}).

The representation of the PACF of FARIMA with $d\in (-1/2,0)$ 
is given by the next theorem.

\begin{thm}\label{thm:2.4}
Let $p, q\in\mathbf{N}\cup\{0\}$ and 
$d\in (-1/2,0)$, and let 
$\{X_n\}$ be a fractional ARIMA $(p,d,q)$ process. 
Then, for $n=1,2,\dots$, 
the representations 
{\rm (\ref{eq:2.6})--(\ref{eq:2.8})} 
of $U_n$, $V_n$ and $\alpha_n$ hold 
with all the sums converging absolutely.
\end{thm}

Here is the result of the type ($d/n$) for FARIMA.

\begin{thm}\label{thm:2.5}
Let $p, q\in\mathbf{N}\cup\{0\}$ and 
$d\in (-1/2, 1/2)\setminus\{0\}$, and let 
$\{X_n\}$ be a fractional ARIMA $(p,d,q)$ process. 
Then the PACF has the asymptotics $(d/n)$.
\end{thm}

This last theorem as well as Theorem \ref{thm:2.3} is an improvement of 
earlier work \cite{I2, I3, IK1} for long-memory or FARIMA processes, 
asserting that
\begin{equation}
\vert \alpha_n\vert\sim \vert d\vert/n \qquad (n\to\infty).
\label{eq:2.23}
\end{equation}
Notice that while the earlier result cannot distinguish 
between the $d$ and $-d$ cases (that is, between positive and negative 
differencing), Theorems \ref{thm:2.3} and \ref{thm:2.5} can. 
In \cite{I2, I3, IK1}, the asymptotic behaviour 
of mean-squared prediction error of the type
\begin{equation}
\delta(n) \sim d^2/n\qquad (n\to\infty)
\label{eq:2.24}
\end{equation}
was first derived and then used in Tauberian arguments 
to prove (\ref{eq:2.23}), 
where
\begin{equation*}
\delta(n):=\frac{\Vert P^{\bot}_{[-n,0]}X_1\Vert^2
-\Vert P^{\bot}_{(-\infty,0]}X_1\Vert^2}
{\Vert P^{\bot}_{(-\infty,0]}X_1\Vert^2}
\qquad(n=1,2,\dots).
\end{equation*}
The proofs of Theorems \ref{thm:2.3} and \ref{thm:2.5}, which are based on 
the representation of the PACF, are more direct and much simpler 
than those of the earlier result.

%%%%%%%%%%%%%%%%%%%%      Section 3.    %%%%%%%%%%%%%%%%%%
\section{Inner products of prediction errors}\label{sec:3}

In this section, we derive some expansions of 
$V_n$ and $U_n$ that we need to prove the representation of the PACF. 
The key is to extend \cite[Theorem 4.1]{I2} properly. 

For $n, k\in\mathbf N$, we define the orthogonal projection operator 
$P_n^k$ by
\begin{equation}
P_n^k:=
\begin{cases}
P_{(-\infty,n-1]} & (k=1,3,5,\dots),\\
P_{[1,\infty)} & (k=2,4,6,\dots).
\end{cases}
\label{eq:3.1}
\end{equation}
It should be noticed that $\{P_n^k : k=1,2,\dots\}$ is merely an 
alternating sequence of 
projection operators, first to the subspace $H_{(-\infty,n-1]}$, then to 
$H_{[1,\infty)}$, and so on.

\begin{thm}\label{thm:3.1}Let $Y_1, Y_2\in H$. 
We assume that $\{X_n\}$ is CND. 
\begin{enumerate}
\item We have
\begin{equation}
\begin{split}
&\left(Y_1, Y_2\right)
=\left((P_1^1)^{\bot}Y_1, (P_1^1)^{\bot}Y_2\right)\\
&\qquad\qquad\qquad +\sum_{k=1}^{\infty}
\left((P_{1}^{k+1})^{\bot}P_{1}^k\cdots P_{1}^1Y_1, 
(P_{1}^{k+1})^{\bot}P_{1}^k\cdots P_{1}^1Y_2\right).
\end{split}
\label{eq:3.2}
\end{equation}
\item We have, for $n=2,3,\dots$,
\begin{equation}
\begin{split}
&\left(P^{\bot}_{[1, n-1]}Y_1, P^{\bot}_{[1, n-1]}Y_2\right)
=\left((P_n^1)^{\bot}Y_1, (P_n^1)^{\bot}Y_2\right)\\
&\qquad\qquad\qquad +\sum_{k=1}^{\infty}
\left((P_{n}^{k+1})^{\bot}P_{n}^k\cdots P_{n}^1Y_1, 
(P_{n}^{k+1})^{\bot}P_{n}^k\cdots P_{n}^1Y_2\right).
\end{split}
\label{eq:3.3}
\end{equation}
\end{enumerate}
\end{thm}

If we put $Y_1=Y_2$ in Theorem \ref{thm:3.1}(2), then it 
reduces to \cite[Theorem 4.1]{I2}.

\begin{proof}
(1)\ The orthogonal decompositions
\begin{align*}
H&=H_{(-\infty,0]}^{\bot}\oplus H_{(-\infty,0]},\\
H&=H_{[1,\infty)}^{\bot}\oplus H_{[1,\infty)}
\end{align*}
of $H$ imply the orthogonal decompositions
\begin{align}
I_H&=
P_{(-\infty,0]}^{\bot}\oplus P_{(-\infty,0]},
\label{eq:3.4}\\
I_H&=
P_{[1,\infty)}^{\bot}\oplus P_{[1,\infty)}
\label{eq:3.5}
\end{align}
of the identity map $I_H$, respectively. 
Repeated use of (\ref{eq:3.4}) and (\ref{eq:3.5}) yields, 
for $m=2,3,\dots$,
\begin{equation*}
\begin{split}
&\left(Y_1, Y_2\right)
=\left((P_1^1)^{\bot}Y_1, (P_1^1)^{\bot}Y_2\right)\\
&\qquad +\sum_{k=1}^{m-1}
\left((P_1^{k+1})^{\bot}P_1^k\cdots P_1^1Y_1, 
(P_1^{k+1})^{\bot}P_1^k\cdots P_1^1Y_2\right) + R^m_1,
\end{split}
\end{equation*}
where $R^m_1:=\left(P_1^m\cdots P_1^1Y_1, P_1^m\cdots P_1^1Y_2\right)$. 
Since (CND) and von Neumann's alternating 
projection theorem (see \cite[Theorem 9.20]{P}) imply
$$
\underset {m\to\infty}{\mbox{s-lim}}\ 
P_1^m\cdots P_1^1=0,
$$
we have $\lim_{m\to\infty}R^m_1=0$, whence (\ref{eq:3.2}). 

\noindent(2)\ For $n=2,3,\dots$, we have the orthogonal decompositions
\begin{align*}
H^{\bot}_{[1,n-1]}
&=H_{(-\infty,n-1]}^{\bot}\oplus 
\left(H_{[1,n-1]}^{\bot}\cap H_{(-\infty,n-1]}\right),\\
H^{\bot}_{[1,n-1]}
&=H_{[1,\infty)}^{\bot}\oplus 
\left(H_{[1,n-1]}^{\bot}\cap H_{[1,\infty)}\right)
\end{align*}
of $H^{\bot}_{[1,n-1]}$, which in turn imply the orthogonal decompositions
\begin{align}
P_{[1,n-1]}^{\bot}&=
P_{(-\infty,n-1]}^{\bot}\oplus P_{[1,n-1]}^{\bot}P_{(-\infty,n-1]},
\label{eq:3.6}\\
P_{[1,n-1]}^{\bot}&=
P_{[1,\infty)}^{\bot}\oplus P_{[1,n-1]}^{\bot}P_{[1,\infty)}
\label{eq:3.7}
\end{align}
of $P^{\bot}_{[1,n-1]}$, respectively. 
Using (\ref{eq:3.6}) and (\ref{eq:3.7}) repeatedly, we find that, 
for $m=2,3,\dots$,
\begin{equation*}
\begin{split}
&\left(P^{\bot}_{[1,n-1]}Y_1, P^{\bot}_{[1,n-1]}Y_2\right)
=\left((P_n^1)^{\bot}Y_1, (P_n^1)^{\bot}Y_2\right)\\
&\qquad\qquad\qquad +\sum_{k=1}^{m-1}
\left((P_n^{k+1})^{\bot}P_n^k\cdots P_n^1Y_1, 
(P_n^{k+1})^{\bot}P_n^k\cdots P_n^1Y_2\right) + R^m_n,
\end{split}
\end{equation*}
where $R^m_n:=(P_{[1,n-1]}^{\bot}P_n^m\cdots P_n^1Y_1, 
P_{[1,n-1]}^{\bot}P_n^m\cdots P_n^1Y_2)$. 
From (IPF) implied by (CND) (see \cite[Theorem 2.3]{IK2}) and 
the alternating projection theorem, we get
$$
\underset {m\to\infty}{\mbox{s-lim}}\ 
P_n^m\cdots P_n^1
=P_{[1,n-1]},
$$
whence
$$
\lim_{m\to\infty}\Vert P_{[1,n-1]}^{\bot}P_n^m\cdots P_n^1Y_i
\Vert
=\Vert P_{[1,n-1]}^{\bot}P_{[1,n-1]}Y_i\Vert=0\qquad(i=1,2).
$$
Thus $\lim_{m\to\infty}R^m_n=0$, so that (\ref{eq:3.3}) follows.
\end{proof}

\begin{remark}
If $\{X_n\}$ is CND, then by the same arguments as above we see that
\[
P_{[1,n-1]}^{\bot}=(P^1_n)^{\bot}+(P^2_n)^{\bot}P^1_n+
(P^3_n)^{\bot}P^2_nP^1_n+\cdots.
\]
\end{remark}

Assuming (PND), we define
\begin{equation*}
b^m_{j}:=\sum_{k=0}^{m}c_{k}a_{j+m-k}
\qquad(m, j\in\mathbf{N}\cup\{0\}).
\end{equation*}
Notice that $b^m_j$ here is equal to that in 
\cite{IK1} but it corresponds to $b^{m+1}_{j-1}$ in \cite{I2, I3}. 

Recall $U_n$ and $V_n$ from \S 2. 
Here are their representations 
in terms of the AR and MA coefficients.

\begin{thm}\label{thm:3.2}
We assume {\rm (PND)} and $\sum_{n=0}^{\infty}\vert a_n\vert<\infty$. 
Then, for $n=1,2,\dots$,
\begin{align}
&U_n=(c_0)^2\sum_{k=1}^{\infty}\sum_{p=0}^{\infty}d_{k}(n,p)d_{k-1}(n,p),
\label{eq:3.8}\\
&V_n=(c_0)^2\sum_{k=0}^{\infty}\sum_{p=0}^{\infty}d_{k}(n,p)^2,
\label{eq:3.9}
\end{align}
where, for $n\in\mathbf{N}$ and $p\in\mathbf{N}\cup\{0\}$, 
$d_0(n,p):=\delta_{p0}$, 
\begin{equation}
d_1(n,p):=\sum_{m_1=0}^{\infty}a_{n+m_1+p}c_{m_1},
\label{eq:3.10}
\end{equation}
and, for $k=2,3,\dots$,
\begin{equation}
d_k(n,p)
:=\sum_{m_{k-1}=0}^{\infty}a_{n+m_{k-1}}
\sum_{m_{k-2}=0}^{\infty}b^{m_{k-1}}_{n+m_{k-2}}
\cdots \sum_{m_{1}=0}^{\infty}b^{m_{2}}_{n+m_{1}}
\sum_{v=0}^{\infty}b^{m_{1}}_{n+p+v}
c_{v},
\label{eq:3.11}
\end{equation}
all the sums converging absolutely.
\end{thm}

\begin{proof}
(Compare the proof of \cite[Theorem 4.5]{I2}.) 
Notice that $\{X_n\}$ satisfies (IPF) 
(see \cite{I2}, Proposition 4.2 and Theorem 3.1), whence (CND) 
(see \cite{IK2}, Theorem 2.3). 
Hence, it follows from Theorem \ref{thm:3.1} that, for $n=1,2,\dots$,
\begin{gather}
\begin{split}
&U_n
=\left((P_n^2)^{\bot}P_n^1X_n, (P_n^2)^{\bot}X_0\right)\\
&\quad\qquad\qquad\qquad+\sum_{k=2}^{\infty}
\left((P_{n}^{k+1})^{\bot}P_{n}^{k}\cdots P_{n}^1X_n, 
(P_{n}^{k+1})^{\bot}P_{n}^{k}\cdots P_{n}^2X_0\right),
\end{split}\label{eq:3.12}\\
V_n
=\Vert (P_n^1)^{\bot}X_n\Vert^2
+\sum_{k=1}^{\infty}
\Vert (P_{n}^{k+1})^{\bot}P_{n}^{k}\cdots P_{n}^1X_n\Vert^2.
\label{eq:3.13}
\end{gather}

Let $n\in\mathbf{N}$. Suppose that $k$ is even and $\ge 2$. 
By \cite[Theorem 4.4]{I2}, we have, for $n=1,2,\dots$ and $m=0,1,\dots$,
\begin{gather*}
P_{(-\infty,n-1]}X_{m+n}=\sum_{j=0}^{\infty}b_{n+j}^{m}X_{-j}
\qquad(\mbox{mod $H_{[1,n-1]}$ if $n\ge 2$}),\\
P_{[1,\infty)}X_{-m}=\sum_{j=0}^{\infty}b_{n+j}^{m}X_{j+n}
\qquad(\mbox{mod $H_{[1,n-1]}$ if $n\ge 2$}),
\end{gather*}
whence
\begin{equation*}
\begin{split}
&P_{n}^{k}\cdots P_{n}^1X_{n}
=c_0\sum_{m_{k-1}=0}^{\infty}a_{n+m_{k-1}}\sum_{m_{k-2}=0}^{\infty}
b^{m_{k-1}}_{n+m_{k-2}}\\
&\qquad\qquad\qquad
\cdots \sum_{m_{1}=0}^{\infty}b^{m_{2}}_{n+m_{1}}
\sum_{m_{0}=0}^{\infty}b^{m_{1}}_{n+m_{0}}
X_{m_0+n}\quad(\mbox{mod $H_{[1,n-1]}$ if $n\ge 2$}).
\end{split}
\end{equation*}
Since we restrict to (PND), $\{X_n\}$ has no deterministic component 
in the Wold decomposition and 
it permits the moving-average representation (MA), 
where the orthonormal system $\{\xi_j: j\in\mathbf{Z}\}$ of $H$ 
satisfies
\begin{equation*}
H_{(-\infty,m]}=H_{(-\infty,m]}(\xi)
\qquad(m\in\mathbf{Z})
\end{equation*}
with $H_{(-\infty,m]}(\xi)$ being the closed subspace of $H$ spanned by 
$\{\xi_{j}: -\infty<j\le m\}$ (see \cite[Chapter II]{Ro}, 
\cite[\S 5.7]{BD}). 
Since
$$
P_{(-\infty,n-1]}^{\bot}X_{m+n}=\sum_{j=0}^{m}c_{m-j}\xi_{j+n}
\qquad(m=0,1,\dots),
$$
we have
\begin{equation*}
\begin{split}
&(P_{n}^{k+1})^{\bot}P_{n}^{k}\cdots P_{n}^1X_{n}
=c_0\sum_{m_{k-1}=0}^{\infty}a_{n+m_{k-1}}\sum_{m_{k-2}=0}^{\infty}
b^{m_{k-1}}_{n+m_{k-2}}\\
&\qquad\qquad\qquad\qquad\qquad\qquad
\cdots \sum_{m_{1}=0}^{\infty}b^{m_{2}}_{n+m_{1}}
\sum_{m_{0}=0}^{\infty}b^{m_{1}}_{n+m_{0}}
\sum_{j=0}^{m_0}c_{m_0-j}\xi_{j+n},
\end{split}
\end{equation*}
so that
\begin{equation*}
\left((P_{n}^{k+1})^{\bot}P_{n}^{k}\cdots P_{n}^1X_{n}, \xi_{p+n}\right)
=
\begin{cases}
c_0d_{k}(n,p)&\qquad(p=0,1,\dots),\\
0&\qquad(p=-1,-2,\dots).
\end{cases}
\end{equation*}
Arguing similarly,
\begin{equation*}
\left(\xi_{p+n}, (P_{n}^{k+1})^{\bot}P_{n}^{k}\cdots P_{n}^2X_{0}\right)
=
\begin{cases}
c_0d_{k-1}(n,p)&\qquad(p=0,1,\dots),\\
0&\qquad(p=-1,-2,\dots).
\end{cases}
\end{equation*}
Thus, from the Parseval equality, we get
\begin{gather}
\begin{split}
&\left((P_{n}^{k+1})^{\bot}P_{n}^{k}\cdots P_{n}^1X_{n}, 
(P_{n}^{k+1})^{\bot}P_{n}^{k}\cdots P_{n}^2X_{0}\right)\\
&\qquad
\qquad\qquad\qquad\qquad=(c_0)^2\sum_{p=0}^{\infty}d_{k}(n,p)d_{k-1}(n,p),
\end{split}
\label{eq:3.14}\\
\Vert (P_{n}^{k+1})^{\bot}P_{n}^{k}\cdots P_{n}^1X_{n}\Vert^2
=(c_0)^2\sum_{p=0}^{\infty}d_{k}(n,p)^2.
\label{eq:3.15}
\end{gather}

Similarly, we have (\ref{eq:3.14}) and (\ref{eq:3.15}) for 
$k$ odd, and also
\begin{gather}
\left((P_{n}^{2})^{\bot}P_{n}^1X_{n}, (P_{n}^{2})^{\bot}X_{0}\right)
=(c_0)^2d_1(n,0)=(c_0)^2\sum_{p=0}^{\infty}d_{1}(n,p)d_{0}(n,p),
\label{eq:3.16}\\
\Vert (P_{n}^{1})^{\bot}X_{n}\Vert^2
=(c_0)^2=(c_0)^2\sum_{p=0}^{\infty}d_{0}(n,p)^2.
\label{eq:3.17}
\end{gather}
The assertions (\ref{eq:3.8}) and (\ref{eq:3.9}) now follow if we 
substitute (\ref{eq:3.14}) and (\ref{eq:3.16}) into (\ref{eq:3.12}), and 
(\ref{eq:3.15}) and (\ref{eq:3.17}) into (\ref{eq:3.13}).
\end{proof}

We write $\sum^{\infty-}$ to indicate that the sum does 
not necessarily converge absolutely, i.e., 
$\sum_{k=m}^{\infty-}:=\lim_{M\to\infty}\sum_{k=m}^{M}$. We need the 
next variant of Theorem \ref{thm:3.2} when we consider the fractional 
ARIMA$(p,d,q)$ processes with $-1/2<d<0$.

\begin{thm}\label{thm:3.3}
We assume {\rm (PND)}. Then 
the representations\/ $(\ref{eq:3.8})$ and\/ $(\ref{eq:3.9})$ 
still hold 
if all the summations $\sum^{\infty}$ in\/ $(\ref{eq:3.10})$ and\/ 
$(\ref{eq:3.11})$ are replaced by $\sum^{\infty-}$ and if\/ 
$\sum_{0}^{\infty}\vert a_n\vert<\infty$ is replaced by the two conditions 
$\sum_{0}^{\infty}\vert c_n\vert<\infty$ and 
$\sum_{0}^{\infty}\vert a_n\vert^2<\infty$.
\end{thm}

\begin{proof}
By \cite[Proposition 4.2 and Theorem 3.1]{I2} and 
\cite[Theorem 2.3]{IK2}, 
the conditions (PND) and $\{a_n\}\in l^2$ 
imply (IPF), whence (CND). 
Moreover, by \cite[Proposition 2.1]{IK1}, we have, for $n\in\mathbf{N}$ 
and $m\in\mathbf{N}\cup\{0\}$,
\begin{gather*}
P_{(-\infty,n-1]}X_{m+n}=\sum_{j=0}^{\infty-}b_{n+j}^{m}X_{-j}
\qquad(\mbox{mod $H_{[1,n-1]}$ if $n\ge 2$}),\\
P_{[1,\infty)}X_{-m}=\sum_{j=0}^{\infty-}b_{n+j}^{m}X_{j+n}
\qquad(\mbox{mod $H_{[1,n-1]}$ if $n\ge 2$}).
\end{gather*}
Using these equalities, we can prove the theorem as in the proof of 
Theorem \ref{thm:3.2}. We omit the details.
\end{proof}

%%%%%%%%%%%%%%%%%%%%      Section 4.    %%%%%%%%%%%%%%%%%%

\section{Proof of Theorem \ref{thm:2.1}}\label{sec:4}

First we give some necessary and sufficient conditions for (BC). 
Recall that the spectral measure of a stationary process is assumed to have 
infinite support.

\begin{thm}\label{thm:4.1}
For a stationary process $\{X_n\}$, the following conditions are 
equivalent:
\begin{enumerate}
\item $\{X_n\}$ is PND and satisfies both 
$\sum_{0}^{\infty}\vert a_n\vert<\infty$ 
and $\sum_{0}^{\infty}\vert c_n\vert<\infty$;
\item $\{X_n\}$ has a positive continuous spectral density and satisfies 
$\sum_{0}^{\infty}\vert c_n\vert<\infty$;
\item $\{X_n\}$ has a positive continuous spectral density and 
satisfies $\sum_{0}^{\infty}\vert a_n\vert<\infty$;
\item $\{X_n\}$ satisfies {\rm (BC)};
\item $\{X_n\}$ has summable PACF: 
$\sum_{0}^{\infty}\vert \alpha_n\vert<\infty$.
\end{enumerate}
\end{thm}

\begin{proof}
Notice that if 
$\{X_n\}$ has a positive continuous spectral density $\Delta(\cdot)$ 
on $[-\pi,\pi]$, 
then it is PND since 
$\int_{-\pi}^{\pi}\vert\log \Delta(\theta)\vert d\theta<\infty$ 
holds. 

Suppose (1). We write $D(e^{i\theta})$ for the nontangential limit of 
$D(z)$, 
i.e., 
\begin{equation*}
D(e^{i\theta})=\lim_{r\to 1-0}D(re^{i\theta})
=\sum\nolimits_{n=0}^{\infty}c_ne^{in\theta}\qquad 
(-\pi\le\theta\le \pi).
\end{equation*}
Since $\{c_n\}\in l^1$ 
implies the continuity of $D(e^{i\theta})$, 
the spectral density $\Delta(\theta)$ is also continuous by the 
equality 
$\Delta(\theta)=2\pi\vert D(e^{i\theta})\vert^2$. 
Letting $r\to 1-0$ in
\begin{equation*}
\left(\sum\nolimits_{n=0}^{\infty}c_nr^ne^{in\theta}\right)
\left(\sum\nolimits_{n=0}^{\infty}a_nr^ne^{in\theta}\right)
=-1\qquad (-\pi\le\theta\le \pi),
\end{equation*}
we obtain
\begin{equation*}
\left(\sum\nolimits_{n=0}^{\infty}c_ne^{in\theta}\right)
\left(\sum\nolimits_{n=0}^{\infty}a_ne^{in\theta}\right)
=-1\qquad (-\pi\le\theta\le \pi).
\end{equation*}
This implies that $D(e^{i\theta})$, whence 
$\Delta(\theta)$, has no zeros on $[-\pi,\pi]$. 
Thus $\Delta(\cdot)$ is positive, 
whence (2) and (3) follow.

Suppose (3). In the same way as above, we have
\begin{equation*}
\frac{1}{\Delta(\theta)}=\frac{1}{2\pi}
\left\vert \sum\nolimits_{n=0}^{\infty}a_ne^{in\theta}\right\vert^2\qquad 
(-\pi\le \theta\le \pi).
\end{equation*}
This implies $\sum\nolimits_{n=0}^{\infty}a_ne^{in\theta}\ne 0$ for every 
$\theta\in [-\pi,\pi]$. By Wiener's theorem for absolutely convergent 
Fourier series (cf.~Lemma 11.6 in \cite{Ru}), 
we obtain $\{c_n\}\in l^1$ 
(cf.~Berk \cite{Berk}, page 493). 
Thus (1) follows. The proof of the implication (2) $\Rightarrow$ (1) 
is similar.

By (\ref{eq:2.1}), (2) implies (4). Conversely, we assume (4). 
Then we have $\{a_n\}\in l^1$, 
whence (3), by 
the arguments in Baxter \cite{Ba2}, pp.~139--140, which involve the 
Wiener--L\'evy theorem.

The equivalence between (4) and (5) is Baxter's theorem (\cite{Ba1,Ba2}; see 
also \cite[Chapter 5]{Si3}). 
 This completes the proof.
\end{proof}

We put
\begin{equation*}
B(n):=\sum_{v=0}^{\infty}\vert c_va_{n+v}\vert\qquad 
(n\in\mathbf{N}\cup\{0\}).
\end{equation*}
For $n,k,u,v\in\mathbf{N}\cup\{0\}$, we define $D_{k}(n,u,v)$ recursively by
\begin{equation*}
\left\{
\begin{aligned}
&D_0(n,u,v):=\delta_{uv},\\
&D_{k+1}(n,u,v):=\sum\nolimits_{w=0}^\infty B(n+v+w)D_k(n,u,w)
\end{aligned}
\right.
\end{equation*}
(see \cite[\S 2.3]{IK2}). We have, for example,
\begin{equation*}
D_3(n,u,v)=\sum_{v_{1}=0}^\infty
\sum_{v_2=0}^\infty B(n+v+v_{1})B(n+v_{1}+v_{2})B(n+v_{2}+u).
\end{equation*}

\begin{prop}\label{prop:4.2}
We assume either {\rm (BC)} or {\rm (O($1/n$))}. 
Then, for $k,n,v\in\mathbf N\cup\{0\}$,
\begin{equation*}
\sum_{u=0}^{\infty}D_{k}(n,u,v)<\infty\quad 
\mbox{and}
\quad
\sum_{u=0}^{\infty}D_{k}(n,u,v)^2<\infty,
\end{equation*}
respectively. In particular, we have 
$D_{k}(n,u,v)<\infty$ for $k,n,u,v\in \mathbf N\cup\{0\}$.
\end{prop}

In view of Theorem \ref{thm:4.1}, we can prove 
Proposition \ref{prop:4.2} 
in the same way as that of \cite[Lemma 2.7]{IK2}, whence we omit it.

Recall $\beta(n)$ from ($\ast$) and $d_k(n,p)$ from Theorem 
\ref{thm:3.2}.

\begin{prop}\label{prop:4.3}
We assume either {\rm (BC)} or {\rm (O($1/n$))}. 
Then we have, for $n\in\mathbf{N}$ and $p\in\mathbf{N}\cup\{0\}$, 
\begin{align}
d_1(n,p)&=\beta(n+p),
\label{eq:4.1}\\
d_2(n,p)&=\sum_{m_1=0}^{\infty}\beta(m_1+n)\beta(m_1+n+p),
\label{eq:4.2}
\end{align}
and, for $k=3,4,\dots$,
\begin{equation}
\begin{split}
&d_k(n,p)
=\sum_{m_{k-1}=0}^{\infty}
\beta(m_{k-1}+n)
\sum_{m_{k-2}=0}^{\infty}\beta(m_{k-1}+m_{k-2}+n) \\
&\qquad\cdots \sum_{m_2=0}^{\infty}\beta(m_3+m_2+n)
\sum_{m_1=0}^{\infty}
\beta(m_2+m_1+n)\beta(m_1+n+p),
\end{split}
\label{eq:4.3}
\end{equation}
the sums converging absolutely.
\end{prop}

\begin{proof}
By Proposition \ref{prop:4.2}, we can use the Fubini--Tonelli theorem to 
exchange the order of sums in (\ref{eq:3.11}), and we get 
(\ref{eq:4.2}) and (\ref{eq:4.3}) as in the proof of 
\cite[Theorem 4.6]{I2}.
\end{proof}

\begin{prop}\label{prop:4.4}
We assume either {\rm (BC)} or {\rm (O($1/n$))}. Then, 
for $i, j\in \mathbf{N}$,
\begin{equation}
\sum_{p=0}^{\infty}d_i(n,p)d_j(n,p)=d_{i+j}(n,0)\qquad(n=1,2,\dots).
\label{eq:4.4}
\end{equation}
\end{prop}

\begin{proof}
For simplicity, we give details for the case $i=j=4$ only. The general case 
can be treated in the same way. 
From Proposition \ref{prop:4.3}, we have, 
for $n=1,2,\dots$ and $p=0,1,\dots$,
\begin{equation}
\begin{split}
&d_4(n,p)
=\sum_{m_{1}=0}^{\infty}\nolimits
\beta(m_{1}+n)
\sum_{m_{2}=0}^{\infty}\beta(m_{1}+m_{2}+n) \\
&\qquad\qquad\quad\sum_{m_3=0}^{\infty}\nolimits
\beta(m_2+m_3+n)\beta(m_3+p+n).
\end{split}
\label{eq:4.5}
\end{equation}
By Proposition \ref{prop:4.2} and the Fubini theorem, 
\begin{equation}
\begin{split}
&d_4(n,p)
=\sum_{m_{3}=0}^{\infty}\nolimits\beta(m_3+p+n)
\sum_{m_{2}=0}^{\infty}\nolimits\beta(m_{2}+m_{3}+n) \\
&\qquad\qquad\quad\sum_{m_1=0}^{\infty}\nolimits
\beta(m_1+m_2+n)\beta(m_{1}+n).
\end{split}
\label{eq:4.6}
\end{equation} 
Writing $(m_7,m_6,m_5)$ for $(m_1,m_2,m_3)$ in (\ref{eq:4.5}), 
we get
\begin{equation}
\begin{split}
&d_4(n,p)
=\sum_{m_{7}=0}^{\infty}\nolimits
\beta(m_{7}+n)
\sum_{m_{6}=0}^{\infty}\nolimits\beta(m_{6}+m_{7}+n) \\
&\qquad\qquad\quad\sum_{m_5=0}^{\infty}\nolimits
\beta(m_5+m_6+n)\beta(p+m_5+n).
\end{split}
\label{eq:4.7}
\end{equation}
From (\ref{eq:4.6}), (\ref{eq:4.7}) and 
the Fubini theorem, 
\begin{equation*}
\begin{split}
&\sum_{p=0}^{\infty}\nolimits d_4(n,p)d_4(n,p)
=\sum_{m_4=0}^{\infty}\nolimits d_4(n,m_4)d_4(n,m_4)\\
&=\sum_{m_4=0}^{\infty}\nolimits\left\{\sum_{m_{7}=0}^{\infty}
\nolimits\beta(m_{7}+n)
\sum_{m_{6}=0}^{\infty}\nolimits\beta(m_{6}+m_{7}+n)\right. \\
&\qquad\qquad\qquad\qquad\qquad\left.\sum_{m_5=0}^{\infty}\nolimits
\beta(m_5+m_6+n)\beta(m_4+m_5+n)\right\}\\
&\qquad\times\left\{\sum_{m_{3}=0}^{\infty}\nolimits\beta(m_3+m_4+n)
\sum_{m_{2}=0}^{\infty}\nolimits\beta(m_{2}+m_{3}+n)\right. \\
&\qquad\qquad\qquad\qquad\qquad\left.\sum_{m_1=0}^{\infty}\nolimits
\beta(m_1+m_2+n)\beta(m_{1}+n)\right\},
\end{split}
\end{equation*}
which is equal to
\begin{equation*}
\begin{split}
&\sum_{m_{7}=0}^{\infty}\nolimits\beta(m_{7}+n)
\sum_{m_{6}=0}^{\infty}\nolimits\beta(m_{6}+m_{7}+n)\sum_{m_5=0}^{\infty}
\beta(m_5+m_6+n) \\
&\quad\sum_{m_4=0}^{\infty}\nolimits\beta(m_4+m_5+n)
\sum_{m_{3}=0}^{\infty}\nolimits\beta(m_3+m_4+n)
\sum_{m_{2}=0}^{\infty}\nolimits\beta(m_{2}+m_{3}+n)\\
&\quad\quad\sum_{m_1=0}^{\infty}\nolimits\beta(m_1+m_2+n)\beta(m_{1}+n)\\
&=d_8(n,0).
\end{split}
\end{equation*}
Thus the desired result for $i=j=4$ follows.
\end{proof}

\noindent {\it Proof of Theorem \ref{thm:2.1}.\/}\ 
By Proposition \ref{prop:4.3}, we see that
\[
d_k(n)=d_k(n,0)\qquad (k,n\in\mathbf{N}).
\]
From this, Proposition \ref{prop:4.4} and Theorem \ref{thm:3.2}, 
the theorem follows. 
\qed

%%%%%%%%%%%%%%%%%%%%      Section 5.    %%%%%%%%%%%%%%%%%%
\section{Proofs of Theorems \ref{thm:2.3}--\ref{thm:2.5}}\label{sec:5}

\begin{prop}\label{prop:5.1}
For $d\in (0,1/2)$ and $\ell\in\mathcal{R}_0$, 
we assume $(L(d,\ell))$. 
\begin{enumerate}
\item It holds that
\begin{equation*}
\beta_n\sim \frac{\sin(\pi d)}{\pi}n^{-1}\qquad(n\to\infty).
\end{equation*}
\item The condition {\rm (O($1/n$))} holds. 
More precisely, we have
\begin{equation*}
\sum_{v=0}^{\infty}\vert c_va_{n+v}\vert
\sim \frac{\sin(\pi d)}{\pi}n^{-1}\qquad(n\to\infty).
\end{equation*}
\item For $s\ge 0$ and $u\ge 0$, it holds that
\begin{equation*}
\beta([ns]+[nu]+n)\sim \frac{\sin(\pi d)}{\pi(s+u+1)}n^{-1}\qquad(n\to\infty).
\end{equation*}
\item For $r\in (1,\infty)$, there exists $N_1\in\mathbf{N}$ such that
\begin{equation}
\vert \beta([ns]+[nu]+n)\vert
\le \frac{r\sin(\pi d)}{\pi(s+u+1)}n^{-1}\quad 
(s\ge 0,\ u\ge 0,\ n\ge N_1).
\label{eq:5.1}
\end{equation}
\end{enumerate}
\end{prop}

\begin{proof}
The assertions (1) and (2) follow from \cite[Proposition 4.3]{I1}. 
Since we have 
$[ns]+[nu]+n\sim n(s+u+1)$ as $n\to\infty$, (3) follows from (1). 
Let $r\in (1,\infty)$. Then $n/([ns]+[nu]+n) \to 1/(s+u+1)$ 
as $n\to\infty$, uniformly in $s\ge 0$ and $u\ge 0$ 
(cf.~\cite[Theorem 1.5.2]{BGT}), so that there exists 
$N_2\in\mathbf{N}$ such that
$$
\frac{1}{([ns]+[nu]+n)}\le \frac{r^{1/2}}{n(s+u+1)}\qquad
(s\ge 0,\ u\ge 0,\ n\ge N_2),
$$
while, from (1), there exists $N_3\in\mathbf{N}$ such that
$$
\vert  \beta_n\vert\le \frac{r^{1/2}\sin(\pi d)}{\pi n}\qquad
(n\ge N_3).
$$
If we put $N_1:=\max(N_2, N_3)$, then, for $s\ge 0$, $u\ge 0$, 
$n\ge N_1$, we have
\begin{equation*}
\vert \beta([ns]+[nu]+n)\vert
\le \frac{r^{1/2}\sin(\pi d)}{\pi([ns]+[nu]+n)}
\le \frac{r\sin(\pi d)}{\pi(s+u+1)}n^{-1}.
\end{equation*}
Thus (4) follows.
\end{proof}

Recall $d_k(n)$ from (\ref{eq:2.3})--(\ref{eq:2.5}). 
For $k=1,2,\dots$, we define the constant $\tau_k$, 
which is equal to $f_k(0)$ in \cite{IK2}, by
\begin{equation*}
\tau_1=\frac{1}{\pi},\quad
\tau_2=\frac{1}{{\pi}^2}\int_{0}^{\infty}\frac{ds_1}{(s_1+1)(s_1+1)}
=\frac{1}{{\pi}^2},
\end{equation*}
and, for $k=3,4,\dots$,
\begin{equation*}
\tau_k=\frac{1}{{\pi}^{k}}\int_0^{\infty}ds_{k-1}
\cdots\int_{0}^{\infty}ds_1\frac{1}{(s_{k-1}+1)}
\left\{\prod_{m=1}^{k-2}\frac{1}{(s_{m+1}+s_{m}+1)}\right\}
\frac{1}{(s_{1}+1)}.
\end{equation*}

\begin{prop}\label{prop:5.2}
For $d\in (0,1/2)$ and $\ell\in\mathcal{R}_0$, 
we assume $(L(d,\ell))$. 
\begin{enumerate}
\item For $r\in (1,\infty)$ and $N_1\in\mathbf{N}$ 
satisfying $(\ref{eq:5.1})$, 
\begin{equation}
\vert d_k(n)\vert
\le n^{-1}\{r\sin(\pi d)\}^k \tau_k\qquad
(u\ge 0,\ k\in\mathbf{N},\ n\ge N_1).
\label{eq:5.2}
\end{equation}
\item For $k\in\mathbf{N}$ and $u\ge 0$,
\begin{equation}
d_k(n)\sim n^{-1}\{\sin(\pi d)\}^k \tau_k\qquad(n\to\infty).
\label{eq:5.3}
\end{equation}
\end{enumerate}
\end{prop}

\begin{proof}
Let $k\ge 3$ and write
\begin{equation*}
\begin{split}
d_k(n)&=\int_0^{\infty}ds_{k-1}
\cdots\int_{0}^{\infty}ds_1\beta([s_{k-1}]+n)\\
&\qquad\times
\left\{\prod_{m=1}^{k-2}\beta([s_{m+1}]+[s_{m}]+n)\right\}
\times \beta([s_{1}]+n)\\
&=n^{k-1}\int_0^{\infty}ds_{k-1}
\cdots\int_{0}^{\infty}ds_1\beta([ns_{k-1}]+n)\\
&\qquad\times
\left\{\prod_{m=1}^{k-2}\beta([ns_{m+1}]+[ns_{m}]+n)\right\}
\times \beta([ns_{1}]+n).
\end{split}
\end{equation*}
Applying Proposition \ref{prop:5.1} 
and the dominated convergence theorem to this, we 
obtain (\ref{eq:5.2}) and (\ref{eq:5.3}). 
The cases $k=1,2$ can be treated in a similar fashion. 
\end{proof}

\begin{prop}\label{prop:5.3}
For $k=1,2,\dots$, we have $\tau_k\le \pi^{-2}$.
\end{prop}

\begin{proof}
Let $T$ be the linear bounded operator on $L^2((0,\infty),du)$ 
defined by
$$
Tg(u):=\int_0^{\infty}\frac{1}{u+v}g(v)dv.
$$
Then, by Hilbert's theorem 
(cf.~\cite[Theorems 316 and 317]{HLP}), the operator norm 
$\Vert T\Vert$ is equal to $\pi$. Hence, for 
the inner product $(\cdot,\cdot)$ of $L^2((0,\infty),du)$ 
and $f(x):=1/(1+x)$, 
we have $\tau_k\le\pi^{-k}(f,T^{k-2}f)\le \pi^{-2}$, 
yielding the proposition.
\end{proof}

\noindent {\it Proof of Theorem \ref{thm:2.3}.\/}\ 
To apply the dominated convergence theorem, we choose 
$r>1$ so that $0<r\sin(\pi d)<1$. 
Then, by Proposition \ref{prop:5.3}, 
$\sum_{k=1}^{\infty}\tau_{2k-1}\{r\sin(\pi d)\}^{2k-1}<\infty$. 
Hence Proposition \ref{prop:5.2} and the dominated convergence 
theorem yield
\begin{equation}
\lim_{n\to\infty}n\sum_{k=1}^{\infty}d_{2k-1}(n)
=\sum_{k=1}^{\infty}\tau_{2k-1}\sin^{2k-1}(\pi d).
\label{eq:5.4}
\end{equation}
By Lemma \ref{lem:5.4} below, the right-hand side is equal to 
$d$. Thus, by Corollary \ref{cor:2.2}, ($d/n$) follows. \qed

\begin{lem}\label{lem:5.4}
For $\vert x\vert<1$, we have 
$\sum_{k=1}^{\infty}\tau_{2k-1}x^{2k-1}=\pi^{-1}\arcsin x$.
\end{lem}

\begin{proof}
(Compare the proof of \cite[Lemma 6.5]{I2}.) 
For $0<d<1/2$, let $\{Y_n:n\in\mathbf{Z}\}$ be a fractional 
ARIMA$(0,d,0)$ process such that $E[(Y_0)^2]=\Gamma(1-2d)/\Gamma^2(1-d)$. 
We denote by $c'_n$, $a'_n$, and $\alpha'_n$ 
the MA and AR coefficients, and PACF of $\{Y_n\}$, 
respectively. Then we have, for $n=0,1,\dots$,
\begin{equation*}
c'_n=\frac{\Gamma(n+d)}{\Gamma(n+1)\Gamma(d)},\qquad
a'_n=\frac{\Gamma(n-d)d}{\Gamma(n+1)\Gamma(1-d)}
\end{equation*}
(see, e.g., \cite[\S 13.2]{BD}). 
We define $d'_k(n)$ similarly. 
Then since $\{Y_n\}$ satisfies (L$(d,\ell')$) with $\ell'\equiv 1/\Gamma(d)$, 
it follows from (\ref{eq:5.4}) that
\begin{equation*}
\lim_{n\to\infty}n\sum_{k=1}^{\infty}d'_{2k-1}(n)
=\sum_{k=1}^{\infty}\tau_{2k-1}\{\sin(\pi d)\}^{2k-1}.
\end{equation*}
However, since $\alpha'_n=d/(n-d)$, Corollary \ref{cor:2.2} gives
\begin{equation*}
\lim_{n\to\infty}n\sum_{k=1}^{\infty}d'_{2k-1}(n)
=\lim_{n\to\infty}n\alpha'_{n}=d.
\end{equation*}
Combining, we obtain $\sum_{1}^{\infty}\tau_{2k-1}\sin^{2k-1}(\pi d)=d$. 
The lemma follows if we substitute $\pi^{-1}\arcsin x$ 
with $0<x<1$ for $d$ and use analytic continuation. \end{proof}

\begin{remark}From Lemma \ref{lem:5.4}, it follows that
\begin{equation*}
\tau_{2k-1}=\frac{1}{\pi}
\cdot\frac{(2k-2)!}{2^{2k-2}((k-1)!)^2(2k-1)}\qquad(k=1,2,\dots).
\end{equation*}
\end{remark}

\ 

\noindent {\it Proof of Theorem \ref{thm:2.4}.\/}\ 
Let $\{X_n\}$ be a fractional ARIMA$(p,d,q)$ 
process with spectral density (\ref{eq:2.14}) with (\ref{eq:2.15}). 
We assume that $-1/2<d<0$. 
Then (\ref{eq:2.16}) and (\ref{eq:2.17}) imply 
$\{c_n\}\in l^1$ and 
$\{a_n\}\in l^2$, 
respectively, so that we can use Theorem \ref{thm:3.3}.
Let $d_k(n,p)$ be as in Theorem \ref{thm:3.3}. 
Recall $\phi_n$, $\psi_n$ and $\beta_{-}(n)$ from \S 2. 
By \cite[Theorem 3.3]{IK1}, we have the same conclusions as those in 
Proposition \ref{prop:4.3} with $\beta(n)$ replaced by $\beta_{-}(n)$. 
Notice that, in \cite[Theorem 3.3]{IK1}, the results 
are stated for $n\ge 2$ but we can prove the case $n=1$ 
in the same way. 
The equality (\ref{eq:4.4}) holds in the same way as the proof of 
Proposition \ref{prop:4.4}. 
Hence the theorem follows from Theorem \ref{thm:3.3}.
\qed

\ 

\noindent {\it Proof of Theorem \ref{thm:2.5}.\/}\ 
If $0<d<1/2$, then ($d/n$) follows immediately from 
Theorem \ref{thm:2.3}. 
We assume that $-1/2<d<0$. Let $q:=d+1\in (1/2,1)$ as in \S 2. 
Then, from (\ref{eq:2.21}), (\ref{eq:2.22}) and \cite[Proposition 4.3]{I1}, 
it follows that 
\[
\beta_-(n)\sim -\frac{\sin(\pi q)}{\pi}n^{-1}\qquad (n\to\infty).
\]
Running through the same arguments as those in the proof of 
Theorem \ref{thm:2.3}, we see that
\begin{equation*}
\lim_{n\to\infty}n\alpha_n
=\lim_{n\to\infty}n\sum_{k=1}^{\infty}d_{2k-1}(n)
=-\sum_{k=1}^{\infty}\tau_{2k-1}\sin^{2k-1}(\pi q)
=d.
\end{equation*}
Thus, again, ($d/n$) holds.
\qed

%%%%%%%%%%%%%%%%%%%%      Section 6.    %%%%%%%%%%%%%%%%%%
\section{Model with regularly varying autocovariance function}\label{sec:6}

In this section, we apply the representation of PACF to a 
stationary process $\{X_n\}$ 
which has regularly varying autocovariance function. 
We will also assume that $\{X_n\}$ 
is PND and satisfies the 
following conditions (cf.~\cite[\S 2]{I2}):
\begin{gather*}
\mbox{\itshape $c_n\ge 0$ for all $n\ge 0$;} \tag{{\rm C1}}\\
\mbox{\itshape $\{c_n\}$ is eventually decreasing to zero;} 
\tag{{\rm C2}}\\
\mbox{\itshape $\{a_n\}$ is eventually decreasing to zero.} 
\tag{{\rm A1}}
\end{gather*}

Notice that (C1) and (A1) imply $\{a_n\}\in l^1$ 
(see \cite[Proposition 4.3]{I2}). 
In \cite{I2}, the extra condition
\begin{equation*}
\mbox{\itshape $\{a_n-a_{n+1}\}$ is eventually decreasing 
to zero} \tag{{\rm A2}}
\end{equation*}
is also required but we do not need it here. 
By \cite[Theorem 7.3]{I2}, 
$\{X_n\}$ satisfies (C1)--(A1) (and also (A2)) if
\begin{equation*}
\begin{split}
&\mbox{there exists a finite Borel measure $\sigma$ on 
$[0,1)$}\\
&\mbox{such that 
$\gamma_n=\int_{0}^{1}t^{\vert n\vert}\sigma(dt)\quad (n\in\mathbf{Z})$.}
\end{split}
\tag{RP}
\end{equation*}
This property is called {\it reflection positivity\/} or 
{\it $T$-positivity}, 
which originates in quantum field theory; 
see, e.g., Osterwalder and Schrader \cite{OS}, Hegerfeldt \cite{He} and 
Okabe \cite{O}. 
A prototype of such a process is $\{X_n\}$ with 
$\gamma_n=(1+\vert n\vert)^{-(1-2d)}$, $-\infty<d<1/2$, 
which we consider in the Example below.

Let $\ell\in\mathcal{R}_0$, and choose a positive constant $B$ so large 
that $\ell(\cdot)$ is locally bounded on $[B,\infty)$ 
(see \cite[Corollary 1.4.2]{BGT}). 
When we say $\int^{\infty}\ell(s)ds/s=\infty$, it means that 
$\int_B^{\infty}\ell(s)ds/s=\infty$. 
If so, then we define another slowly varying function $\tilde{\ell}$ by
\begin{equation}
\tilde{\ell}(x):=\int_B^{x}\frac{\ell(s)}{s}ds\qquad(x\ge B)
\label{eq:1.25}
\end{equation}
(see \cite[\S 1.5.6]{BGT}). 
The asymptotic behaviour of $\tilde{\ell}(x)$ as $x\to\infty$ does not 
depend on the choice of $B$ since we have assumed that 
$\int^{\infty}\ell(s)ds/s=\infty$. 

Here is the result on the asymptotic behaviour of $\alpha_n$.

\begin{thm}\label{thm:6.1}Let $-\infty<d<1/2$ and 
$\ell\in\mathcal{R}_0$. We assume {\rm (PND)}, 
{\rm (C1)}, {\rm (C2)}, {\rm (A1)}, 
and
\begin{equation}
\gamma_n\sim n^{2d-1}\ell(n)\qquad(n\to\infty).
\label{eq:6.2}
\end{equation}
\begin{enumerate}
\item If $0<d<1/2$, then $(d/n)$ holds;
\item if $d=0$ and $\int^{\infty}\ell(s)ds/s=\infty$, then
\begin{equation}
\alpha_n
\sim n^{-1}\frac{\ell(n)}{2\tilde{\ell}(n)}
\qquad(n\to\infty);
\label{eq:6.3}
\end{equation}
\item if $d=0$ with $\int^{\infty}\ell(s)ds/s<\infty$ or $-\infty<d<0$, 
then
\begin{equation}
\alpha_n
\sim \frac{n^{2d-1}\ell(n)}
{\sum_{-\infty}^{\infty}\gamma_k}
\qquad(n\to\infty).
\label{eq:6.4}
\end{equation}
\end{enumerate}
\end{thm}

This theorem is an improvement of \cite[Theorem 2.1]{I2} where 
only $\vert\alpha_n\vert$ is considered 
with additional assumption (A2).

We need some propositions to prove the theorem above.

\begin{prop}\label{prop:6.2}
Let $\ell\in\mathcal{R}_0$. If 
$\int^{\infty}\ell(s)ds/s=\infty$, then 
$\ell(n)/\tilde{\ell}(n)$ tends to $0$ as $n\to\infty$. If 
$\int^{\infty}\ell(s)ds/s<\infty$, then 
$\ell(n)$ tends to $0$ as $n\to\infty$.
\end{prop}

Proposition \ref{prop:6.2} follows immediately from 
\cite[Proposition 1.5.9a]{BGT}.

\begin{prop}\label{prop:6.3}Let $\ell\in\mathcal{R}_0$ and 
$-\infty<d\le 0$. We assume {\rm (PND)}, 
{\rm (C1)}, {\rm (C2)}, {\rm (A1)}, and $(\ref{eq:6.2})$.
\begin{enumerate}
\item If $d=0$ and $\int^{\infty}\ell(s)ds/s=\infty$, 
then
\begin{gather}
c_n\sim n^{-1}\ell(n)\{2\tilde{\ell}(n)\}^{-1/2}\qquad(n\to\infty),
\label{eq:6.5}\\
a_n\sim n^{-1}\ell(n)\{2\tilde{\ell}(n)\}^{-3/2}\qquad(n\to\infty),
\label{eq:6.6}\\
\beta_n\sim n^{-1}\ell(n)\{2\tilde{\ell}(n)\}^{-1}\qquad(n\to\infty).
\label{eq:6.7}
\end{gather}
\item If $d=0$ with $\int^{\infty}\ell(s)ds/s<\infty$ or $-\infty<d<0$, 
then
\begin{gather}
c_n\sim n^{2d-1}\ell(n)
\left\{\sum\nolimits_{-\infty}^{\infty}\gamma_k\right\}^{-1/2}
\qquad(n\to\infty),
\label{eq:6.8}\\
a_n\sim n^{2d-1}\ell(n)
\left\{\sum\nolimits_{-\infty}^{\infty}\gamma_k\right\}^{-3/2}
\qquad(n\to\infty),
\label{eq:6.9}\\
\beta_n\sim 
n^{2d-1}\ell(n)\left\{\sum\nolimits_{-\infty}^{\infty}\gamma_k\right\}^{-1}
\qquad(n\to\infty).
\label{eq:6.10}
\end{gather}
\end{enumerate}
\end{prop}

\begin{proof}
The assertions (\ref{eq:6.5}) and (\ref{eq:6.6}) follow from 
\cite[Theorem 5.2]{I2}. Using them, we obtain (\ref{eq:6.7}) 
(see \cite[(6.19)]{I2}). 
The assertions (\ref{eq:6.8}) and (\ref{eq:6.9}) follow from 
\cite[Theorem 5.3]{I2}. From them, we get (\ref{eq:6.10}) 
(see the proof of \cite[Theorem 6.7]{I2}). 
\end{proof}

\begin{prop}\label{prop:6.4}
Let $\ell\in\mathcal{R}_0$ and $-\infty<d\le 0$. 
We assume {\rm (PND)}, 
{\rm (C1)}, {\rm (C2)}, {\rm (A1)}, and $(\ref{eq:6.2})$.
\begin{enumerate}
\item For every $R\in (1,\infty)$, there exists 
$N\in\mathbf{N}$ such that
\begin{equation}
\left\vert \frac{\beta([ns]+[nu]+n)}{\beta(n)}\right\vert
\le \frac{R}{(s+u+1)}\qquad
(s\ge 0,\ u\ge 0,\ n\ge N).
\label{eq:6.11}
\end{equation}
\item For every $r\in (0,1)$, there exists 
$N\in\mathbf{N}$ such that
\begin{equation}
\vert \beta([ns]+[nu]+n)\vert
\le \frac{r}{\pi(s+u+1)}n^{-1}\qquad
(s\ge 0,\ u\ge 0,\ n\ge N).
\label{eq:6.12}
\end{equation}
\end{enumerate}
\end{prop}

\begin{proof}By Proposition \ref{prop:6.3} and 
\cite[Theorem 1.5.2]{BGT}, we have
$$
\beta([ns]+[nu]+n)/\beta(n)\to (s+u+1)^{2d-1}\qquad(n\to\infty)
$$
uniformly in $s\ge 0$ and $u\ge 0$. Since
$$
(s+u+1)^{2d-1}\le (s+u+1)^{-1},
$$
(1) follows. By Propositions \ref{prop:6.2} and \ref{prop:6.3}, we have
\begin{equation}
\lim_{n\to\infty}n\beta_n=0.
\label{eq:6.13}
\end{equation}
This and (1) show (2). 
\end{proof}

Recall $d_k(n)$ from (\ref{eq:2.3})--(\ref{eq:2.5}) 
and $\tau_k$ from \S 5.

\begin{prop}\label{prop:6.5}
Let $\ell\in\mathcal{R}_0$ and $-\infty<d\le 0$. 
We assume {\rm (PND)}, 
{\rm (C1)}, {\rm (C2)}, {\rm (A1)}, and $(\ref{eq:6.2})$.
\begin{enumerate}
\item Let $r\in (0,1)$ and $R\in (1,\infty)$. 
Choose $N\in\mathbf{N}$ so that both 
$(\ref{eq:6.11})$ and $(\ref{eq:6.12})$ hold. Then 
$\vert d_k(n)/\beta_n\vert
\le \pi Rr^{k-1}\tau_k$ for 
$k\in\mathbf{N},\ n\ge N$.
\item For $k\ge 2$, 
$\lim_{n\to\infty} d_k(n)/\beta_n=0$.
\end{enumerate}
\end{prop}

\begin{proof}
In the same way as the proof of (\ref{eq:6.13}), we see that 
(O($1/n$)) holds. 
We assume $k\ge 3$; the cases $k=1,2$ can be treated in a similar way. 
We write
\begin{equation*}
\begin{split}
\frac{d_k(n)}{\beta(n)}&=\int_0^{\infty}ds_{k-1}
\cdots\int_{0}^{\infty}ds_1\frac{\beta([ns_{k-1}]+n)}{\beta(n)}\\
&\times
\left\{\prod_{m=1}^{k-2}n\beta([ns_{m+1}]+[ns_{m}]+n)\right\}
\times n\beta([ns_{1}]+n).
\end{split}
\end{equation*}
From this as well as (\ref{eq:6.11}) and 
(\ref{eq:6.12}), we get (1). 
The assertion (2) follows from (\ref{eq:6.11})--(\ref{eq:6.13}) and 
the dominated convergence theorem.
\end{proof}

\noindent {\it Proof of Theorem \ref{thm:6.1}.\/}\ 
Assume that $0<d<1/2$. Then, from \cite[Theorem 5.1]{I2}, 
(L($d, \ell'$)) holds with 
$\ell'(n):=[\ell(n)/B(d,1-2d)]^{1/2}$. 
Hence ($d/n$) follows 
immediately from Theorem \ref{thm:2.3}.

Next, we assume $-\infty<d\le 0$. 
By (A1) and (C1), $\sum_{v=0}^{\infty}\vert c_va_{v+n}\vert =\beta_n$ 
if $n$ is large enough. 
Hence from (\ref{eq:6.7}), (\ref{eq:6.10}) and Proposition \ref{prop:6.2}, 
(O($1/n$)) holds. 
By Proposition \ref{prop:6.5} and the dominated convergence theorem,
we have
$$
\lim_{n\to\infty}\frac{\sum_{k=1}^{\infty}d_{2k-1}(n)}{\beta_n}
=1+\lim_{n\to\infty}
\sum_{k=2}^{\infty}\frac{d_{2k-1}(n)}{\beta_n}=1.
$$
Therefore, by Corollary \ref{cor:2.2}, we see (\ref{eq:2.10}). 
Thus (2) and (3) follow from (\ref{eq:6.7}) and 
(\ref{eq:6.10}), respectively.
\qed

\begin{exmp}
Let $-\infty<d<1/2$, and let $\{X_n\}$ be a stationary process with 
autocovariance function of the form 
$\gamma_n=(1+\vert n\vert)^{-(1-2d)}$. Then 
$\{X_n\}$ satisfies (RP) (cf.~\cite[Example in \S 7]{I2}). 
Let $\{\alpha_n\}$ be the PACF of $\{X_n\}$. Applying Theorem \ref{thm:6.1} to 
$\{X_n\}$, we get the following result:
\begin{enumerate}
\item if $0<d<1/2$, then we have ($d/n$).
\item if $d=0$, then 
$\alpha_n\sim 1/(2n\log n)$ as $n\to\infty$.
\item if $-\infty<d<0$, then 
$\alpha_n\sim n^{2d-1}\cdot [2\zeta(1-2d)-1]^{-1}$ as 
$n\to\infty$. 
\end{enumerate}
Here $\zeta(s)$ is the Riemann zeta function. 
\end{exmp}
\ 

\noindent{\it Remarks.}\quad 
1.\ 
Recall by (2.3) that $\beta_n$ is the first term on the right 
of (\ref{eq:2.9}). 
By the arguments above, we see that, 
for the processes treated in Theorem \ref{thm:6.1}, 
\begin{equation}
\lim_{n\to\infty}\frac{\alpha_n}{\beta_n}
=
\begin{cases}
\frac{\pi d}{\sin(\pi d)} & (0<d<1/2),\\
1 & (-\infty<d\le 0).
\label{eq:??}
\end{cases}
\end{equation}
We raise, and leave open, the question of how generally this happens.

\noindent 2.\ We note that in Theorems \ref{thm:2.3}, \ref{thm:2.5} and 
\ref{thm:6.1} we have also
\begin{equation*}
\alpha_n\sim \frac{\gamma_n}{\sum_{k=-n}^{n}\gamma_k}\qquad
(n\to\infty)
\end{equation*}
(see \cite[\S 2 and \S 6]{I2} and \cite[\S 5]{I3} for the proofs). 
It would be interesting to know how generally this relation holds.

\

%%%%%%%%%%%%%%%%%%%%      Section 7.    %%%%%%%%%%%%%%%%%%
\section{ARMA processes}\label{sec:7}

In this section, we consider the fractional ARIMA$(p,0,q)$ processes, 
that is, the ARMA$(p,q)$ processes. 
Let $p, q\in\mathbf{N}\cup\{0\}$, and let 
$\Phi(z)$ and $\Theta(z)$ be polynomials with real coefficients of degrees 
$p$, $q$, respectively, satisfying (\ref{eq:2.15}). 
Let $\{X_n\}$ be an ARMA$(p,q)$ process with spectral 
density
\begin{equation*}
\Delta(\theta)
=\frac{1}{2\pi}
\frac{\vert \Theta(e^{i\theta})\vert^2}{\vert \Phi(e^{i\theta})\vert^2}
\qquad(-\pi<\theta<\pi).
\end{equation*}
We put $R:=0$ if $q=0$ and 
\begin{equation*}
R:=
\max\left(1/\vert u_1\vert,\dots,1/\vert u_q\vert\right)
\quad \mbox{if $q\ge 1$},
\end{equation*}
where $u_1,\dots,u_q$ are the (complex) zeros of $\Theta(z)$:
\begin{equation*}
\Theta(z)=\mbox{const.}\times (z-u_1)\cdots(z-u_q).
\end{equation*}
From the assumption (\ref{eq:2.15}), we see that $\vert u_k\vert >1$ for 
$k=1,\dots,q$, whence $R\in [0,1)$.

Let $\{\alpha_n\}$ be the 
PACF of the ARMA$(p,q)$ process $\{X_n\}$. 
The next theorem implies that $\alpha_n$ decays exponentially 
as $n\to\infty$.

\begin{thm}\label{thm:7.1}
For every $r>R$, we have
\begin{equation}
\alpha_n=O(r^n)\qquad(n\to\infty).
\label{eq:7.1}\end{equation}
In particular, $\alpha_n$ decays exponentially fast as $n\to\infty$.
\end{thm}

\begin{proof}
The Szeg\"o function $D(z)$ of $\{X_n\}$ is given by 
$D(z)=\Theta(z)/\Phi(z)$ for $\vert z\vert<1$. 
Hence $-\Phi(z)/\Theta(z)=\sum_{0}^{\infty}a_nz^n$, so that, for every $r>R$,
\begin{equation*}
a_n=O(r^n)\qquad(n\to\infty).
\end{equation*}
Treating $D(z)=\Theta(z)/\Phi(z)$ similarly, 
we see that $c_n$ also decays exponentially as $n\to\infty$, 
in particular, $\{c_n\}\in l^1$. 
Therefore, we see that, for every $r>R$,
\[
\beta_n=O(r^n)\qquad(n\to\infty).
\]
The assertion (\ref{eq:7.1}) follows easily from this and 
Theorem \ref{thm:2.1}. 
Since $r$ can be chosen so that $R<r<1$, (\ref{eq:7.1}) 
implies exponential decay of $\alpha_n$.
\end{proof}

\section*{Acknowledgements}
The author expresses his gratitude to 
Professor Nick Bingham and Dr.\ Yukio Kasahara.  Their 
comments led to substantial improvements of the paper, and the proposed new definition of short memory based on Baxter's theorem is due to them. 
He is also grateful to three referees for their 
constructive comments.

%%%%%%%%%%%%%%%%%        References       %%%%%%%%%%%%%%%%%%%%

 \end{document}